\documentclass[10pt,reqno]{amsart}
\usepackage{amscd,amsmath,amsopn,amssymb,amsthm,multicol}
\usepackage{hyperref}
\DeclareMathOperator{\Ad}{Ad}

\DeclareMathOperator{\C}{constant}
\DeclareMathOperator{\Ric}{Ric}
\DeclareMathOperator{\g}{\mathfrak{g}}
\DeclareMathOperator{\hh}{\mathfrak{h}}
\DeclareMathOperator{\pp}{\mathfrak{p}}

\newtheorem{theorem}{Theorem}
\newtheorem{lemma}{Lemma}

\newtheorem{prop}{Proposition}

\begin{document}

\title[Einstein invariant metrics...]
{Invariant Einstein metrics on some homogeneous spaces of classical
Lie groups}

\author{Andreas Arvanitoyeorgos, V.V. Dzhepko and Yu.G.~Nikonorov}

\address{University of Patras, Department of Mathematics, GR-26500 Rion, Greece}
\email{arvanito@math.upatras.gr}
\address{Rubtsovsk industrial institute, ul. Traktornaya, 2/6, Rubtsovsk, 658207, Russia}
\email{J\_Valera\_V@mail.ru}
\address{Rubtsovsk industrial institute, ul. Traktornaya, 2/6, Rubtsovsk, 658207, Russia}
\email{nik@inst.rubtsovsk.ru}

\begin{abstract}
A Riemannian manifold $(M,\rho)$ is called Einstein if the metric
$\rho$ satisfies the condition $\Ric (\rho)=c\cdot \rho$ for some
constant $c$.   This paper is devoted to the investigation of
$G$-invariant Einstein metrics with additional symmetries,
on some homogeneous spaces $G/H$ of classical groups.
As a consequence, we obtain new invariant Einstein metrics on some
Stiefel manifolds $SO(n)/SO(l)$, and on the symplectic analogues $Sp(n)/Sp(l)$.
Furthermore, we show that for any positive integer $p$ there exists a Stiefel manifold $SO(n)/SO(l)$ and
a homogenous space $Sp(n)/Sp(l)$ which admit at least $p$
$SO(n)$ (resp. $Sp(n)$)-invariant Einstein metrics.

\medskip
\noindent 2000 Mathematical Subject Classification: 53C25, 53C30.

\medskip
\noindent Key words and phrases: Riemannian manifolds, Homogeneous
spaces, Einstein metrics, Stiefel manifolds.
\end{abstract}


\maketitle

\section*{Introduction}

A Riemannian manifold $(M,\rho)$ is called Einstein if the metric
$\rho$ satisfies the condition $\Ric (\rho)=c\cdot \rho$ for some
real constant $c$. A detailed exposition on Einstein manifolds can
be found in the book of A.~Besse \cite{Bes}, and more recent results
on homogeneous Einstein manifolds can be found in the survey of
M.~Wang \cite{Wang2}. General existence results are hard to obtain.
Among the first important attempts are the works of G.~Jensen
\cite{Je1} and M.~Wang, W.~Ziller \cite{WZ1}.  Recently, a new
existence approach was introduced by C.~B\"ohm, M.~Wang, and W.~Ziller
\cite{BWZ}, \cite{Bo}. The above existence results were used by
C.~B\"ohm and M.~Kerr in \cite{BK} to show that every compact simply
connected homogeneous space up to dimension 11 admits at least one
invariant Einstein metric. It is known (\cite{Bes}, \cite{BK},
\cite{WZ1}) that in dimension 12 there are examples of compact
simply connected homogeneous spaces that do not admit any invariant
Einstein metrics.

The structure of the set of invariant Einstein metrics on a given
homogeneous space is still not very well understood in general.  The
situation is only clear for few classes of homogeneous spaces, such
as isotropy irreducible homogeneous spaces, low dimensional
examples, certain flag manifolds, and some other special types of
homogeneous spaces (\cite{Ar}, \cite{Bes}, \cite{Ki},
\cite{LNF}, \cite{Sak}). For an arbitrary compact homogeneous space
$G/H$ it is not clear if the set of invariant Einstein metrics (up
to isometry and up to scaling) is finite or not (cf. \cite{WZ2}).  A
finiteness conjecture states that this set is in fact finite if the
isotropy representation of $G/H$ consists of pairwise inequivalent
irreducible components (\cite[p. 683]{BWZ}).

Let $G$ be a compact Lie group and $H$ a closed subgroup so that $G$
acts almost effectively on $G/H$.
 In this paper we investigate
$G$-invariant metrics on $G/H$ with additional symmetries. More
precisely, let $K$ be a closed subgroup of $G$ with $H\subset
K\subset G$, and suppose that $K=L'\times H'$, where
$\{e_{L'}\}\times H'=H$. It is clear
that $K\subset N_G(H)$, the normalizer of
$H$ in $G$.  If we denote $L=L'\times\{e_{H'}\}$, then
the group $\widetilde G=G\times L$ acts on $G/H$ by
$(a,b)\cdot gH=agb^{-1}H$, and the isotropy subgroup at $eH$ is $\widetilde
H=\{(a,b): ab^{-1}\in H\}$.

Later on it will be shown that the set $\mathcal M^{\widetilde G}$ of
$\widetilde G$-invariant metrics on $\widetilde G/\widetilde H$ is a
subset of $\mathcal M^G$, the set of $G$-invariant metrics on $G/H$. Therefore,
it would be simpler to search for invariant Einstein metrics on
$\mathcal M^{\widetilde G}$. In this way we obtain existence results
for Einstein metrics for certain quotients.

We apply this method for the case of Stiefel manifolds
$SO(n)/SO(n-k)$.  Note that the simplest case
$S^{n-1}=SO(n)/SO(n-1)$ is an irreducible symmetric space, therefore
it admits up to scale a unique invariant Einstein metric. Concerning
history, it was S.~Kobayashi \cite{Ko} who proved first the
existence of an invariant Einstein metric on $T_1S^n=SO(n)/SO(n-2)$.  Later on, A.
Sagle \cite{Sa} proved that the Stiefel manifolds $SO(n)/SO(n-k)$
admit at least one homogeneous invariant Einstein metric. For $k\ge
3$ G.~Jensen \cite{Je2} found a second metric.  In the same work he
also proved that $Sp(n)/Sp(n-k)$ admits at least two homogeneous
invariant Einstein metrics.  Einstein metrics on $SO(n)/SO(n-2)$ are
completely classified. If $n=3$ the group $SO(3)$ has a unique
Einstein metric.  If $n\ge 5$ it was shown by A.~Back and W.Y.~Hsiang
\cite{BH} that $SO(n)/SO(n-2)$ admits exactly one homogeneous
invariant Einstein metric.  The same result was obtained by M. Kerr
\cite{Kerr1}. The Stiefel manifold $SO(4)/SO(2)$ admits exactly two
invariant Einstein metrics which follows from the classification of
$5$-dimensional homogeneous Einstein manifolds due to
D.V.~Alekseevsky, I.~Dotti, and C.~Ferraris \cite{ADF}. We also refer to
\cite[p.727-728]{BWZ} for further discussion. For $k\ge 3$ there is
no obstruction for existence of more than two homogeneous invariant
Einstein metrics on Stiefel manifolds $SO(n)/SO(n-k)$.

In particular we prove the following:

\begin{theorem}\label{T1} If $s>1$ and $l>k\ge 3$ then the Stiefel manifold $SO(sk+l)/SO(l)$ admits at least four $SO(sk+l)\times (SO(k))^s$-invariant
Einstein metrics, two of which are Jensen's metrics.
\end{theorem}

Analogously we prove:
\begin{theorem}\label{T2} If $s>1$ and $l\ge k\ge 1$ then the space $Sp(sk+l)/Sp(l)$ admits at least four $Sp(sk+l)\times (Sp(k))^s$-invariant Einstein metrics, two of which are Jensen's metrics.
\end{theorem}

We also prove the following:

\begin{theorem}\label{T3} For any positive integer $p$ there exists a Stiefel manifold $SO(n)/SO(l)$ and
a homogenous space $Sp(n)/Sp(l)$ which admit at least $p$
$SO(n)$ (resp. $Sp(n)$)-invariant Einstein metrics.
\end{theorem}

We remark that in fact there are other homogeneous spaces for which the
number of invariant Einstein metrics can be at least a prescribed
number. Indeed, if a compact homogeneous space $G/H$ admits two
distinct invariant Einstein metrics, then the product of $m$
copies of this space admits at least
$m+1$ distinct Einstein metrics invariant under the natural action of $G^m$.
There are analogous examples in the class of non-product homogeneous spaces.
For instance, in \cite[p.62]{DZ} it is shown in particular that
the groups $SU(2n)$, $SU(2n+3)$, $Sp(2n)$, and $Sp(2n+1)$ admit
at least $n+1$ distinct left-invariant Einstein metrics, whereas
the groups $SO(2n)$ and $SO(2n+1)$
admit at least $3n-2$ left-invariant Einstein metrics.

We note also that the methods of this paper (after some minor revisions) can be used
for obtaining new Einstein invariant metric on homogeneous spaces of the groups
$SU(n)$. It is interesting to note that the case of unitary groups $SU(n)$
is more tractable (in some sense) than the cases of orthogonal groups $SO(n)$ and
symplectic groups $Sp(n)$. For instance, one can compare the Einstein equations
and its solutions for the spaces $SU(k_1+k_2+k_3)/S(U(k_1)\times U(k_2) \times U(k_3))$
(see \cite{Ki} and \cite{Ar}) with the Einstein equations
and its solutions for the spaces
$SO(k_1+k_2+k_3)/SO(k_1)\times SO(k_2)\times SO(k_3)$
and
$Sp(k_1+k_2+k_3)/Sp(k_1)\times Sp(k_2)\times Sp(k_3)$
\cite{LNF}.

\medskip
The paper is organized as follows:  The basic construction for searching for
invariant Einstein metrics with additional symmetries on $G/H$ is presented in Section 1, where
we also clarify the meaning of such symmetries.  Metrics with this property are described
for some homogeneous spaces of the groups $SO(n)$ and $Sp(n)$ (Section 2, Lemma \ref{L2}).
In Section 3 we compute the scalar curvature for these metrics (Proposition \ref{main}), and the
variational approach to the Einstein metrics is given in Proposition \ref{variat1}.
 In Section \ref{SectionJen}, as an application
of our construction, we obtain Jensen's invariant Einstein metrics on the Stiefel manifold $SO(k_1+k_2)/SO(k_2)$, and on $Sp(k_1+k_2)/Sp(k_2)$.
In Section \ref{Section4} we investigate $SO(k_1+k_2+k_3)\times SO(k_1)\times SO(k_2)$-invariant Einstein metrics
on $SO(k_1+k_2+k_3)/SO(k_3)$ (resp. $Sp(k_1+k_2+k_3)\times Sp(k_1)\times Sp(k_2)$-invariant Einstein metrics
on $Sp(k_1+k_2+k_3)/Sp(k_3)$).  Another construction for searching for invariant Einstein metrics on $SO(sk+l)/SO(l)$ and
$Sp(sk+l)/Sp(l)$ is given in Section \ref{Section5}.  Finally, in Section \ref{Section6} the proofs of the main results
are given.

\medskip
The second and the third authors are supported in part by the Russian Foundation for~Basic
Research(grant N 05-01-00611-a) and by
the~State Maintenance Program for the~Leading Scientific Schools  of the~Russian
Federation (grants NSH-8526.2006.1) The third author is supported
by the~State Maintenance Program for
Young Russian Scientists of the~Russian
Federation (grant MD-5179.2006.1)

The third author is grateful for the hospitality of the Department
of Mathematics of the University of Patras-Greece, where the final
version of this paper has been prepared.

\section{The main construction}\label{Section1}

Let $G$ be a compact Lie group and $H$ a closed subgroup so that $G$
acts almost effectively on $G/H$.
Let $\g$, $\hh$ be the Lie algebras of $G$ and $H$, and let
$\g=\hh\oplus\pp$ be a reductive decomposition of $\g$ with respect
to some $\Ad(G)$-invariant inner product of $\g$. The orthogonal
complement $\pp$ can be identified with the tangent space
$T_{eH}G/H$.  Any $G$-invariant metric $\rho$ of $G/H$ corresponds
to an $\Ad(H)$-invariant inner product $(\cdot ,\cdot )$
on $\pp$ and vice-versa.  For $G$ semisimple, the negative of the Killing form $B$ of $\g$ is an
$\Ad (G)$-invariant inner product on $\g$, therefore we can choose the above
decomposition with respect to this form.
We will use such a decomposition later on.
Moreover, the restriction $\langle \cdot, \cdot\rangle =-B|_{\pp}$ is
an $\Ad (G)$-invariant inner product on $\pp$, which generates a $G$-invariant metric
on $G/H$ called {\it standard}.

The normalizer $N_G(H)$ of $H$ in $G$  acts on $G/H$ by $(a,
gH)\mapsto ga^{-1}H$. For a fixed $a$ this action induces a
$G$-equivariant diffeomorphism $\varphi _a: G/H\to G/H$. Note that
if $a\in H$ this diffeomorphism is trivial, so the action of the
gauge group $N_G(H)/H$ is well defined. However, it is simpler from
technical point of view to use the action of $N_G(H)$. Let $\rho$ be
a $G$-invariant metric of $G/H$ with corresponding inner product
$(\cdot ,\cdot )$. Then the diffeomorphism $\varphi _a$ is an
isometry of $(G/H, \rho)$ if and only if the operator $\Ad
(a)|_{\pp}$ is orthogonal with respect to $(\cdot ,\cdot )$.

Let $K$ be a closed subgroup of $G$ with $H\subset
K\subset G$ such that $K=L'\times H'$, where
$\{e_{L'}\}\times H'=H$, and consider $L=L'\times\{e_{H'}\}$. It is clear that
$K\subset N_G(H)$.
  The group $\widetilde G=G\times L$ acts on $G/H$ by
$(a,b)\cdot gH=agb^{-1}H$, and the isotropy at $eH$ is
given as follows:

\begin{lemma}\label{action}
The isotropy subgroup $\widetilde H$ is isomorphic to $K$.
\end{lemma}
\begin{proof}  It is clear that $\widetilde H=\{(a,b)\in G\times L: ab^{-1}\in H\}$.
Let $i: K\hookrightarrow G$ be the inclusion of $K$ in $G$.
Then $i(\{e_{L'}\}\times H')=H$ and $i(L'\times\{e_{H'}\})=L$. Let
$(a,b)\in G\times L$ be such that $ab^{-1}=h\in H$. Then $a=hb$, so
$$(a,b)=(hb, b)=(i(b', e_{H'})i(e_{L'}, h'), i(b', e_{H'}))=
((b', h'), b')\in K\times L'=L'\times H'\times L'.
$$
Thus $\widetilde H$ is identified with a subgroup of $L'\times
H'\times L'$, and it is then obvious that $\widetilde H$ is
isomorphic to $L'\times H'=K$.

\end{proof}

The set $\mathcal M^G$ of $G$-invariant metrics on $G/H$ is finite dimensional.
We consider the subset $\mathcal M^{G, K}$ of  $\mathcal M^G$
  corresponding to $\Ad (K)$-invariant inner products on $\pp$ (and not only $\Ad (H)$-invariant).

  Let $\rho\in\mathcal M^{G, K}$ and $a\in K$. The above diffeomorphism $\varphi _a$
 is an isometry of $(G/H, \rho)$.  The action $\widetilde G$  on $(G/H, \rho)$ is isometric,
 so any metric form $\mathcal M^{G, K}$ can be identified a metric in
 $\mathcal M^{\widetilde G}$ and vice-versa.
 Therefore, we may think of $\mathcal M^{\widetilde G}$ as $\mathcal M^{G, K}$, which is a subset of $\mathcal M^G$.

  Since metrics in $\mathcal M^{G,K}$ correspond to
$\Ad(K)$-invariant inner products on $\pp$, we call these metrics
{\it $\Ad (K)$-invariant metrics on $G/H$}.

 The aim of this work is to apply
 the above construction for $G=SO(n)$ and $ Sp(n)$, and
 prove existence of Einstein metrics in the set $\mathcal M^{G, K}$ for
 various choices of the subgroup $K=L'\times H'$.

 \medskip

Let $n\in \mathbb{N}$ and  $k_1$, $k_2,\dots ,k_s,k_{s+1},\dots
,k_{s+t}$ be natural numbers such that $k_1+\dots +k_s=l$,
$k_{s+1}+\dots +k_{s+t}=m$, $l+m=n$. Let $G=SO(n)$ and $K=L'\times
H'$,  where $L'=SO(k_1)\times \cdots \times SO(k_s)$ and
$H'=SO(k_{s+1})\times \cdots \times SO(k_{t+s})$.  The embedding of
$K$ in $G$ is the standard one.
Analogously, we consider $G=Sp(n)$ and $K=L'\times H'$, where
$L'=Sp(k_1)\times \cdots \times Sp(k_s)$ and $H'=Sp(k_{s+1})\times
\cdots \times Sp(k_{t+s})$.

We note that for $s=0$ we obtain a flag manifold of orthogonal or
symplectic groups. Invariant Einstein metrics on the spaces
$SO(k_1+k_2+k_3)/SO(k_1)\times SO(k_2) \times SO(k_3)$ and
$Sp(k_1+k_2+k_3)/Sp(k_1)\times Sp(k_2) \times Sp(k_3)$, were studied
in \cite{LNF}, \cite{N}.

\section{$\Ad(K)$-invariant metrics on the space $G/H$}

Let $\pp _i$ be the subalgebra $so(k_i)$ (resp. $sp(k_i)$) in $\g$,
$1\leq i\leq s+t$. We note that for $1\leq 1 \leq s$  the submodule $\pp _i$
of $\pp$
is $\Ad(K)$-invariant and $\Ad(K)$-irreducible submodule.
For $1\leq i <j \leq s+t$ we denote by $\pp _{(i,j)}$ the  $\Ad(K)$-invariant
and $\Ad(K)$-irreducible submodule of $\pp$ which is determined by the equality
$so(k_i+k_j)=so(k_i)\oplus so(k_j) \oplus \pp_ {(i,j)}$
(resp. $sp(k_i+k_j)=sp(k_i)\oplus sp(k_j) \oplus \pp _{(i,j)}$),
where $\pp _{(i,j)}$ is orthogonal to
$so(k_i)\oplus so(k_j)$
(resp. $sp(k_i)\oplus sp(k_j)$)
with respect to the Killing form $B$.

Denote by $d_i$ and $d_{(i,j)}$ the dimensions of the modules $\pp _i$
and $\pp _{(i,j)}$ respectively. It is easy to obtain that
$d_i=\frac{k_i(k_i-1)}{2}$, $d_{(i,j)}=k_ik_j$ in the orthogonal
case, and $d_i=2k_i^2+k_i$, $d_{(i,j)}=4k_ik_j$ in the symplectic
case.

In both cases we have a decomposition of $\pp$ into a sum of
$\Ad(K)$-invariant and
$\Ad(K)$-irreducible submodules:
\begin{equation}\label{decomp1}
\pp=\bigoplus\limits_{i=1}^s \pp_ i \bigoplus
\bigoplus\limits_{1\leq i <j\leq s+t} \pp_ {(i,j)}.
\end{equation}

\begin{lemma}\label{L1}
If in the orthogonal case for all
$1\leq i \leq s+t$ we have $k_i\geq 2$, and
there is at most one
$1\leq i\leq s$ such that $k_i=2$,
then there are no pairwise $\Ad(K)$-isomorphic submodules among
$\pp_ i$ ($i=1,\dots,s$) and $\pp_ {(i,j)}$ ($1\leq i <j \leq s+t$).

In the symplectic case there are no pairwise $\Ad(K)$-isomorphic
submodules among $\pp_ i$ ($i=1,\dots,s$) and $\pp_ {(i,j)}$ ($1\leq
i <j \leq s+t$).
\end{lemma}

\begin{proof}
It is clear that $\pp_ i$ and $\pp_ j$ act nontrivially on $\pp_
{(i,j)}$, and that each $\pp_ i$ acts nontrivially on itself. Moreover,
the last action is trivial if and only if $k_i=2$ in the orthogonal
case. Therefore, there are no pairwise $\Ad(K)$-isomorphic
submodules.
\end{proof}

If the assumptions of Lemma \ref{L1} are satisfied, then we have a
complete description of the $\Ad (K)$-invariant metrics on $G/H$.

Let $\rho$ be any $\Ad(K)$-invariant metric on $G/H$ with corresponding
$\Ad(K)$-invariant inner product
$(\cdot , \cdot )$ on $\pp$.

\begin{lemma} \label{L2}
If there are no pairwise
$\Ad(K)$-isomorphic submodules among
$p_i$ and $p_{(i,j)}$,
then
\begin{equation}\label{1}
(\cdot,\cdot)=
\sum\limits_{i=1}^s x_{i}\cdot \langle \cdot,\cdot \rangle |_{\pp _{i}}+
\sum\limits_{1\leq i<j\leq s+t} x_{(i,j)}
\cdot \langle \cdot,\cdot \rangle |_{\pp _{(i,j)}}
\end{equation}
for  positive constants $x_{i}>0$ and $x_{(i,j)}>0$, where
$\langle \cdot,\cdot \rangle =-B|_{\pp}$. Therefore, the
set of $\Ad(K)$-invariant metrics on $G/H$ depends on
$(s+t)(s+t-1)/2+s$ parameters.
\end{lemma}

In the case  of pairwise $\Ad(K)$-isomorphic modules $\pp _{\alpha}$
and $\pp _{\beta}$ the set of $\Ad(K)$-invariant metrics have a more
complicated structure \cite{WZ1}.

\section{The scalar curvature and the Einstein condition}\label{Section3}

Let $\{e^{j}_{\alpha}\}$ be an orthonormal basis of $\pp_ {\alpha}$
with respect to $\langle \cdot,\cdot \rangle$, where $1\leq j \leq
d_{\alpha}$ (here $\alpha$ means any of the symbols of type $i$ or
$(k,i)$). We define the numbers (cf. \cite{WZ1}) $[\alpha \beta
\gamma]$ by the equation
$$
[\alpha \beta \gamma]=
\sum\limits_{i,j,k}
{\langle {\left[
{e^i_\alpha ,e^j_\beta } \right] ,e^k_\gamma } \rangle}^2 ,
$$
where
$i,j,k$ vary from 1 to
$d_{\alpha},d_{\beta},d_{\gamma}$ respectively. The symbols
$[\alpha \beta \gamma]$
are symmetric with respect to all three indices, as follows from the
$\Ad (G)$-invariance of $\langle \cdot,\cdot \rangle$.

For any Lie algebra $\mathfrak q$ we shall use the symbol $B_{\mathfrak q}$ for the Killing form
of $\mathfrak q$.
If a simple algebra $\mathfrak q$ is a subalgebra of a Lie algebra
$\mathfrak r$, then we denote by
${\alpha}_{\mathfrak r}^{\mathfrak q}$ a real number which satisfies the equality
$B_{\mathfrak q}={\alpha}_{\mathfrak r}^{\mathfrak q}\cdot B_{\mathfrak r}|_{\mathfrak q}$.

\begin{lemma} \label{L3}
Let $\mathfrak q\subset \mathfrak r$ be arbitrary subalgebras in
$\g$ with $\mathfrak q$ simple. Consider in $\mathfrak q$ an
orthonormal (with respect to $-B_{\mathfrak r}$) basis $\{ f_j\}$
($1\leq j \leq \dim(\mathfrak q)$). Then
$$
\sum\limits_{j,k=1}^{\dim(\mathfrak q)}(-B_{\mathfrak r}([f_i,f_j],f_k))^2={\alpha}_{\mathfrak r}^{\mathfrak q},
\quad i=1,\dots,\dim(\mathfrak q),
$$
$$
\sum\limits_{i,j,k=1}^{\dim(\mathfrak q)}(-B_{\mathfrak r}([f_i,f_j],f_k))^2={\alpha}_{\mathfrak r}^{\mathfrak q}\cdot
\dim(\mathfrak q),
$$
where ${\alpha}_{\mathfrak r}^{\mathfrak q}$ is determined by the equation
$B_{\mathfrak q}={\alpha}_{\mathfrak r}^{\mathfrak q}\cdot B_{\mathfrak r}|_{\mathfrak q}$.
\end{lemma}

\begin{proof}
By direct computations we have
$$
\sum\limits_{j,k}(-B_{\mathfrak r}([f_i,f_j],f_k))^2=
\sum\limits_{j} -B_{\mathfrak r}([f_i,f_j],[f_i,f_j])=
\frac{1}{{\alpha}_{\mathfrak r}^{\mathfrak q}} \sum\limits_{j}
 -B_{\mathfrak q}([f_i,f_j],[f_i,f_j]).
$$
Let $\widetilde{f}_j=\frac{1}{\sqrt{{\alpha}_{\mathfrak r}^{\mathfrak q}}}f_j$. Then
the vectors $\widetilde{f}_j$ form an orthonormal basis in $\mathfrak q$
with respect to $-B_{\mathfrak q}$ and for $i=1,\dots,\dim(\mathfrak q)$ we have that
$$
1= -B_{\mathfrak q}(\widetilde{f}_i,\widetilde{f}_i)=
\sum\limits_{j}
B_{\mathfrak q}([\widetilde{f}_i,[\widetilde{f}_i ,\widetilde{f}_j]], \widetilde{f}_j)=
$$
$$
=\sum\limits_{j}
-B_{\mathfrak q}([\widetilde{f}_i ,\widetilde{f}_j], [\widetilde{f}_i,\widetilde{f}_j])=
\frac{1}{({\alpha}_{\mathfrak r}^{\mathfrak q})^2}
\sum\limits_{j}
-B_{\mathfrak q}([f_i ,f_j], [f_i,f_j]).
$$
Here we use the definition and the properties of the Killing form.
The first statement of the lemma can be easily obtained from the last two formulas.
The second statement is a direct corollary of the first one.
\end{proof}

Using this lemma we obtain an explicit expression for
$[\alpha \beta \gamma]$.
It is clear that the only non-zero
symbols (up to permutation of indices)
are
$$
[aaa],\quad [a(a,b)(a,b)], \quad [b(a,b)(a,b)],
$$
where $1\leq a <b\leq s+t$, and
$[(a,b)(b,c)(a,b)]$,
with $1\leq a<b<c\leq s+t$.

\begin{lemma}\label{L4}
For the orthogonal case the following relations hold:
$$
[aaa]=
\frac{k_a(k_a-1)(k_a-2)}{2(n-2)},\quad
[a(a,b)(a,b)]=
\frac{k_ak_b(k_a-1)}{2(n-2)},
$$
$$
[b(a,b)(a,b)]=
\frac{k_ak_b(k_b-1)}{2(n-2)},\quad
[(a,b)(b,c)(a,c)]=
\frac{k_ak_bk_c}{2(n-2)}.
$$

For the symplectic case the following relations hold:
$$
[aaa]=\frac{k_a(k_a+1)(2k_a+1)}{n+1},\quad
[a(a,b)(a,b)]=\frac{k_ak_b(2k_a+1)}{n+1},
$$
$$
[b(a,b)(a,b)]=
\frac{k_ak_b(2k_b+1)}{n+1},\quad
[(a,b)(b,c)(a,c)]=
\frac{2k_ak_bk_c}{n+1}.
$$
\end{lemma}

\begin{proof}
We give the proof for the orthogonal case.
For the standard embedding
$so(k)\subset so(n)$
we have
${\alpha}_{so(n)}^{so(k)}=\frac{k-2}{n-2}$
(see e.g. \cite{DZ}).

The first equality $[aaa]= \frac{k_a(k_a-1)(k_a-2)}{2(n-2)}$ follows
from Lemma \ref{L3}. In fact, $d_a=\dim(so(k_a))=k_a(k_a-1)/2$ and
${\alpha}_{so(n)}^{so(k_a)}=\frac{k_a-2}{n-2}$.

To prove the second equality we consider the subalgebra
$so(k_a+k_b)\subset so(n)$. It is clear that $[\pp_ a,\pp_ b]=0$,
$[\pp_ a,\pp_ {(a,b)}]\subset \pp_ {(a,b)}$. According to Lemma
\ref{L3} we have that
$$
[aaa]+[a(a,b)(a,b)]=
\dim(\pp_ a)\cdot
{\alpha}_{so(n)}^{so(k_a+k_b)}=\frac{k_a(k_a-1)(k_a+k_b-2)}{2(n-2)},
$$
which proves the second equality.
The third equality can be obtained analogously.

To prove the forth equality we consider
the subalgebra $so(k_a+k_b+k_c)\subset so(n)$.
It is clear that
$$
\dim(\pp _{(a,b)})\cdot
{\alpha}_{so(n)}^{so(k_a+k_b+k_c)}=
2\left(
[(a,b)a(a,b)]+[(a,b)b(a,b)]+[(a,b)(b,c)(a,c)]
\right),
$$
from which we obtain the last equality.

Similar computations apply for the symplectic case. We only need to
note that for the standard embedding $sp(k)\subset sp(n)$ we have
${\alpha}_{sp(n)}^{sp(k)}=\frac{k+1}{n+1}$.
\end{proof}

According to \cite{WZ1}, the scalar curvature
$S$ of $(\cdot ,\cdot)$ is given by

$$
S((\cdot , \cdot ))=
\frac{1}{2} \sum\limits_{\alpha} \frac{d_{\alpha}}{x_{\alpha}}-
\frac{1}{4}
\sum\limits_{\alpha,\beta,\gamma}
[\alpha \beta \gamma]
\frac{x_{\gamma}}{x_{\alpha}x_{\beta}},
$$
where $\alpha$, $\beta$ and $\gamma$ are arbitrary symbols
of the type $i$ ($1\leq i \leq s$)
or of the type $(i,j)$ ($1\leq i <j\leq s+t$).

For the metric (\ref{1}) this formula takes the following form, for
the orthogonal and for the symplectic case.

\begin{prop}\label{main}
The scalar curvature $S$ of an $\Ad(K)$-invariant metric (\ref{1})
has the form
$$
S=\sum\limits_{a=1}^s\frac{k_a(k_a-1)(k_a-2)}{8(n-2)}\cdot \frac{1}{x_a}+
\frac{1}{2}
\sum\limits_{1\leq a<b\leq s+t} \frac{k_ak_b}{x_{(a,b)}}
$$
$$
-\frac{1}{8(n-2)}
\sum\limits_{1\leq a\leq s, \,a+1\leq b\leq s+t} k_ak_b(k_a-1)
\frac{x_a}{x_{(a,b)}^2}-
\frac{1}{8(n-2)}
\sum\limits_{1\leq a<b\leq s} k_ak_b(k_b-1)
\frac{x_b}{x_{(a,b)}^2}
$$
\begin{equation}\label{2}
-\frac{1}{4(n-2)}
\sum\limits_{1\leq a <b<c \leq s+t} k_ak_bk_c
\left(
\frac{x_{(a,b)}}{x_{(a,c)}x_{(b,c)}}+
\frac{x_{(a,c)}}{x_{(a,b)}x_{(b,c)}}+
\frac{x_{(b,c)}}{x_{(a,b)}x_{(a,c)}}
\right)
\end{equation}
in the orthogonal case, and
$$
S=\sum\limits_{a=1}^s\frac{k_a(k_a+1)(2k_a+1)}{4(n+1)}\cdot \frac{1}{x_a}+
2
\sum\limits_{1\leq a<b\leq s+t} \frac{k_ak_b}{x_{(a,b)}}
$$
$$
-\frac{1}{4(n+1)}
\sum\limits_{1\leq a\leq s,\, a+1\leq b\leq s+t} k_ak_b(2k_a+1)
\frac{x_a}{x_{(a,b)}^2}-
\frac{1}{4(n+1)}
\sum\limits_{1\leq a<b \leq s} k_ak_b(2k_b+1)
\frac{x_b}{x_{(a,b)}^2}
$$
\begin{equation}\label{3}
-\frac{1}{n+1}
\sum\limits_{1\leq a <b<c \leq s+t} k_ak_bk_c
\left(
\frac{x_{(a,b)}}{x_{(a,c)}x_{(b,c)}}+
\frac{x_{(a,c)}}{x_{(a,b)}x_{(b,c)}}+
\frac{x_{(b,c)}}{x_{(a,b)}x_{(a,c)}}
\right)
\end{equation}
in the symplectic case.
\end{prop}

Denote by $\mathcal{M}_1^G$ the set of all $G$-invariant metrics with a fixed
volume element on the space $G/H$.
The following variational principle for
invariant Einstein metrics is well known.

\begin{prop}[\cite{Bes}]\label{variat}
Let $G/H$ be a homogeneous space, where $G$ and $H$
are compact.
Then the $G$-invariant Einstein metrics on the homogeneous space $G/H$
are precisely the critical points
of the scalar curvature functional
$S$
restricted to
$\mathcal{M}_1^G$.
\end{prop}

For the general construction as described in Section 1, the
above variational principle implies the following:

\begin{prop}\label{variat1}
Let $\mathcal M_1^{G,K}$ be the subset of $\mathcal M^{G,K}$ with fixed volume element.
Then a metric in $\mathcal M_1^{G,K}$ is Einstein if and only if it is a critical
point of the scalar curvature functional $S$ restricted to $\mathcal M_1^{G,K}$.
\end{prop}
\begin{proof}
The set $\mathcal M_1^{G,K}$ is precisely the set of $\widetilde G$-invariant metrics with fixed
volume element on $\widetilde G/\widetilde H$.
\end{proof}

The volume condition for the metric (\ref{1}) takes the form

\begin{equation}\label{4}
\prod\limits_{i=1}^sx_i^{d_i}\cdot \prod\limits_{1\leq i<j\leq s+t}
x_{(i,j)}^{d_{(i,j)}}=\C.
\end{equation}

By using Proposition \ref{variat1}
the problem of searching for $\Ad(K)$-invariant Einstein metrics on $G/H$
reduces to a Lagrange-type
problem for the scalar curvature functional $S$
under the  restriction (\ref{4}).

\section{Jensen's metrics}\label{SectionJen}

As a first simple illustration of Proposition \ref{variat1} we will
show that Jensen's metrics (\cite{Je2}) on the Stiefel manifold
$SO(k_1+k_2)/SO(k_2)$ ($k_1\ge 2$) can be obtained. We apply
Proposition \ref{main}, formula (\ref{2}) for $s=1$ and $t=1$.  Then
the scalar curvature reduces to
$$
S=\frac{k_1(k_1-1)(k_1-2)}{8(n-2)}\frac{1}{x_1}+\frac12\frac{k_1k_2}{x_{12}}-\frac{1}{8(n-2)}k_1k_2(k_1-1)
\frac{x_1}{x_{12}{}^2}.
$$
The volume condition (\ref{4}) is $V=x_1^{d_1}x_{12}^{d_{12}}=\C$. By
use of Lagrange method we obtain the equation
$$
(k_1-2)x_{12}^2-2(k_1+k_2-2)x_1x_{12}+(k_2+k_1-1)x_1^2=0.
$$
If $k_1=2$ the above equation has a unique solution $x_{12}=\frac{k_2+1}{2k_2} x_1$.
If $k_1>2$ the equation has two solutions
$$
x_{12}=\frac{k_1+k_2-2\pm \sqrt{(k_1+k_2-2)^2-(k_1-2)(k_1+k_2-1)}}{k_1-2}x_1.
$$
These solutions are $SO(k_1+k_2)\times SO(k_1)$-invariant Einstein metrics on $SO(k_1+k_2)/SO(k_2)$, and were
found by G. Jensen in \cite{Je2}.

\medskip
Now apply Proposition \ref{main}, formula (\ref{3}) for $s=1$ and
$t=1$ in the symplectic case, and perform similar computations we obtain the equation
$$
2(k_1+1)x_{12}^2+2k_2x_1^2-4(k_1+k_2+1)x_1x_{12}+(2k_1+1)x_1^2=0,
$$
which has two solutions
$$
x_{12}=\frac{2(k_1+k_2+1)\pm\sqrt{2}\sqrt{1+k_1+2k_2+2k_1k_2+2k_2^2}}{2(k_1+1)}x_1.
$$
These solutions are $Sp(k_1+k_2)\times Sp(k_1)$-invariant Einstein metrics on $Sp(k_1+k_2)/Sp(k_2)$, which
were also found in \cite{Je2}.
Note that if  $k_1=1$ the above solutions simplify to $x_{12}=\frac12 x_1$ and $x_{12}=\frac{2k_2+3}{2} x_1$, which are
Einstein metrics on the sphere $S^{4k_2+3}$ (cf. \cite[Example p. 612-613]{Je2}).

\section{New examples of Einstein metrics}\label{Section4}

In this section we will investigate $SO(k_1+k_2+k_3)\times
SO(k_1)\times SO(k_2)$-invariant (resp. $Sp(k_1+k_2+k_3)\times
Sp(k_1)\times Sp(k_2)$-invariant) Einstein metrics on the spaces
$SO(k_1+k_2+k_3)/SO(k_3)$ and $Sp(k_1+k_2+k_3)/Sp(k_3)$.  Here $L'=SO(k_1)\times SO(k_2)$
(resp. $Sp(k_1)\times Sp(k_2)$).
By Lemma \ref{L2} these metrics depend on $5$ parameters.

We apply Proposition \ref{main} for $s=2$ and $t=1$ for the
orthogonal case  $SO(k_1+k_2+k_3)/SO(k_3)$ ($k_2>2$) and by the
Lagrange method we obtain the following system of equations

\begin{equation}\label{example1}
\begin{array}{l}
x_2x_{23}^2\Bigl((k_1-2)x_{12}^2x_{13}^2+k_2x_1^2x_{13}^2+k_3x_1^2x_{12}^2\Bigr)=
x_1x_{13}^2\Bigl((k_2-2)x_{12}^2x_{23}^2+k_3x_2^2x_{12}^2+k_1x_2^2x_{23}^2\Bigr),\\
\\
x_{13}\Bigl((k_2-2)x_{12}^2x_{23}^2+k_3x_2^2x_{12}^2+k_1x_2^2x_{23}^2\Bigr)=
x_2x_{23}\Bigl(2(k_1+k_2+k_3-2)x_{12}x_{13}x_{23}-\\
(k_1-1)x_1x_{13}x_{23}-
(k_2-1)x_2x_{13}x_{23}+k_3x_{12}^3-k_3x_{12}x_{13}^2-k_3x_{12}x_{23}^2\Bigr),\\
\\
x_{13}\Bigl(2(k_1+k_2+k_3-2)x_{12}x_{13}x_{23}-(k_1-1)x_1x_{13}x_{23}-(k_2-1)x_2x_{13}x_{23}+k_3x_{12}^3-\\
k_3x_{12}x_{13}^2-k_3x_{12}x_{23}^2\Bigr)=\\
x_{12}\Bigl(2(k_1+k_2+k_3-2)x_{12}x_{13}x_{23}-(k_1-1)x_1x_{12}x_{23}+k_2x_{13}^3-k_2x_{12}^2x_{13}-k_2x_{13}x_{23}^2\Bigr),\\
\\
x_{23}\Bigl(2(k_1+k_2+k_3-2)x_{12}x_{13}x_{23}-(k_1-1)x_1x_{12}x_{23}+k_2x_{13}^3-
k_2x_{12}^2x_{13}-k_2x_{13}x_{23}^2\Bigr)=\\
x_{13}\Bigl(2(k_1+k_2+k_3-2)x_{12}x_{13}x_{23}-(k_2-1)x_2x_{12}x_{13}+k_1x_{23}^3-
k_1x_{12}^2x_{23}-k_1x_{13}^2x_{23}\Bigr).
\end{array}
\end{equation}

If  $x_{13}=x_{23}=z$ then system (\ref{example1})
reduces to the following:

\begin{center}
\begin{equation}
\label{example22}
\begin{array}{lcr}
x_{12}z^2\Bigl((k_1-k_2)x_{12}+(k_2-1)x_2-(k_1-1)x_1\Bigr)=0,\\
\\
z^2\Bigl(\Bigl(((k_1-2)x_2-(k_2-2)x_1)x_{12}^2+(k_2x_1-k_1x_2)x_1x_2\Bigr)z^2+
k_3(x_1-x_2)x_1x_2x_{12}^2\Bigr)=0,\\
\\
z\Bigl(\Bigl((k_2-2)x_{12}^2+(k_1+k_2-1)x_2^2-2(k_1+k_2-2)x_2x_{12}+(k_1-1)x_1x_2\Bigr)z^2+\\
k_3(x_2-x_{12})x_2x_{12}^2\Bigr)=0,\\
\\
z\Bigl(\Bigl(2(k_1+k_2-2)x_{12}-(k_2-1)x_2-(k_1-1)x_1\Bigr)z^2+\\
\Bigl((k_2+k_3)x_{12}-2(k_3+k_2+k_1-2)z+(k_1-1)x_1\Bigr)x_{12}^2\Bigr)=0.
\end{array}
\end{equation}
\end{center}

From the first equation of system (\ref{example22}) we obtain
$x_1=\frac{(k_1-k_2)x_{12}+(k_2-1)x_2}{k_1-1}$, and substituting to the two other equations
 we obtain

\begin{center}
\begin{equation}
\label{example23}
\begin{array}{lcr}
z(x_2-x_{12})\Bigl(((2k_2+k_1-2)x_2-(k_2-2)x_{12})z^2+k_3x_2x_{12}^2\Bigr)=0,\\
\\
z\Bigl(\Bigl((k_1+3k_2-4)x_{12}-2(k_2-1)x_2\Bigr)z^2+\\
+\Bigl((k_3+k_1)x_{12}-2(k_1+k_2+k_3-2)z+(k_2-1)x_2\Bigr)x_{12}^2\Bigr)=0,\\
\\
(k_1-k_2)z^2(x_{12}-x_2)\Bigl(\Bigl((1-k_1)(k_2-2)x_{12}^2+(k_1(k_1-2)+k_2(k_1-k_2)+1)x_{12}x_2+\\
+(k_1(k_2-1)+k_2(k_2-2)+1)x_2^2\Bigr)z^2+k_3\Bigl((k_1-k_2)x_{12}+(k_2-1)x_2\Bigr)x_2x_{12}^2\Bigr)=0.
\end{array}
\end{equation}
\end{center}

If $x_2=x_{12}$, we obtain the two Jensen's solutions found in
Section \ref{Section4}. Suppose that $x_2\neq x_{12}$. Then
from the first equation of system (\ref{example23}) we get
$x_2=\frac{(k_2-2)x_{12}z^2}{k_3x_{12}^2+(2k_2+k_1-2)z^2}$. Substituting
  to (\ref{example23}) we obtain the following
system:

\begin{center}
\begin{equation}
\label{example2}
\begin{array}{lcr}
z^4x_{12}^3(k_2-1)(k_2-2)(k_1-k_2)\Bigl((k_1+k_2-1)z^2+k_3x_{12}^2\Bigr)
\Bigl((k_1+k_2)z^2+k_3x_{12}^2\Bigr)^2=0,\\
\\
(k_3^2+k_1k_3)x_{12}^4+(4k_3-2k_2k_3-2k_1k_3-2k_3^2)zx_{12}^3+(k_1^2-6k_3+5k_2k_3+\\
2k_1k_3-2k_1+k_2^2+2k_1k_2-3k_2+2)z^2x_{12}^2+(8k_1-6k_1k_2-4k_2k_3+12k_2-2k_1^2+\\
4k_3-4k_2^2-2k_1k_3-8)z^3x_{12}+(5k_1k_2-8k_2+k_1^2+4k_2^2-6k_1+4)z^4=0.
\end{array}
\end{equation}
\end{center}

From the first equation of (\ref{example2}) we see that
a solution exists only when $k_1=k_2$.

Let $k_1=k_2=k$ ($k\geq 3$), $k_3=l$, $x_1=x_2=x$, $x_{12}=1$,
$x_{13}=x_{23}=z$. Then the original system (\ref{example1})
reduces to the following:

\begin{center}
\begin{equation}
\label{example3}
\begin{array}{lcr}
z(x-1)\Bigl(((3k-2)x-k+2)z^2+lx \Bigr)=0,\\
z\Bigl((2kx-2x-4k+4)z^2+(4k+2l-4)z-(k-1)x-k-l\Bigr)=0.
\end{array}
\end{equation}
\end{center}

If $x=1$ we obtain Jensen's solutions, so assume
that $x\neq 1$. Then from the first equation of
(\ref{example3}) we find that
$$
x=\frac{(k-2)z^2}{(3k-2)z^2+l},
$$
which implies in particular, that $0<x<1$ for any positive $z$.
Substituting this expression to the second equation of
(\ref{example3}), the Einstein equation reduces to $F_{SO}(z)=0$, where
\begin{center}
\begin{equation}
\label{example4}
\begin{array}{lcr}
F_{SO}(z)=2(5k^2-7k+2)z^4-2(6k^2+3kl-10k-2l+4)z^3+\\
(4k^2+7kl-5k-6l+2)z^2-2l(2k+l-2)z+l(k+l).
\end{array}
\end{equation}
\end{center}

If the equation $F_{SO}(z)=0$ has a positive solution, then we obtain a  new
$SO(2k+l)\times SO(k)\times SO(k)$-invariant Einstein metric
on $SO(2k+l)/SO(l)$. The numbers of new Einstein metrics
for some values of $k,l$ are shown
in Table 1.  However, we can show that there exists an infinite series of new homogeneous
Einstein manifolds, as the next proposition shows.

\begin{table}[t]
\begin{center}
\begin{tabular}{|c||c|c|c|c|c|c|c|c|c|c|c|c|c|c|c|c|c|c|c|}
\hline
$k$&3&4&5&6&7&8&9&10&11&12&13&14&15&16&17&18&19&20\\
$l$& & & & & & & & & & & & & & & & & & \\
\hline \hline
1 &0&0&0&0&0&0&0&0&0&0&0&0&0&0&0&0&0&0\\
\hline
2 &0&0&0&0&0&0&0&0&0&0&0&0&0&0&0&0&0&0\\
\hline
3 &2&0&0&0&0&0&0&0&0&0&0&0&0&0&0&0&0&0\\
\hline
4 &2&2&2&0&0&0&0&0&0&0&0&0&0&0&0&0&0&0\\
\hline
5 &2&2&2&2&2&0&0&0&0&0&0&0&0&0&0&0&0&0\\
\hline
6 &2&2&2&2&2&2&2&0&0&0&0&0&0&0&0&0&0&0\\
\hline
7 &2&2&2&2&2&2&2&2&2&2&0&0&0&0&0&0&0&0\\
\hline
8 &2&2&2&2&2&2&2&2&2&2&2&2&0&0&0&0&0&0\\
\hline
9 &2&2&2&2&2&2&2&2&2&2&2&2&2&2&0&0&0&0\\
\hline
10&2&2&2&2&2&2&2&2&2&2&2&2&2&2&2&2&2&0\\
\hline
11&2&2&2&2&2&2&2&2&2&2&2&2&2&2&2&2&2&2\\
\hline
12&2&2&2&2&2&2&2&2&2&2&2&2&2&2&2&2&2&2\\
\hline
13&2&2&2&2&2&2&2&2&2&2&2&2&2&2&2&2&2&2\\
\hline
14&2&2&2&2&2&2&2&2&2&2&2&2&2&2&2&2&2&2\\
\hline
15&2&2&2&2&2&2&2&2&2&2&2&2&2&2&2&2&2&2\\
\hline
16&2&2&2&2&2&2&2&2&2&2&2&2&2&2&2&2&2&2\\
\hline
17&2&2&2&2&2&2&2&2&2&2&2&2&2&2&2&2&2&2\\
\hline
18&2&2&2&2&2&2&2&2&2&2&2&2&2&2&2&2&2&2\\
\hline
19&2&2&2&2&2&2&2&2&2&2&2&2&2&2&2&2&2&2\\
\hline
20&2&2&2&2&2&2&2&2&2&2&2&2&2&2&2&2&2&2\\
\hline
\end{tabular}

\vspace {5mm} {\bf Table 1.} The number of positive solutions of the equation
$F_{SO}(z)=0$ for various $(k,l)$.  (New $SO(2k+l)\times SO(k)\times
SO(k)$-invariant Einstein metrics on $SO(2k+l)/SO(l)$).
\end{center}
\end{table}

\begin{prop}\label{simplSO}
If $l>k\ge 3$ then the Stiefel manifold $SO(2k+l)/SO(l)$ admits at
least four
 $SO(2k+l)\times SO(k)\times SO(k)$-invariant
Einstein metrics.
\end{prop}
\begin{proof}
We consider the polynomial (\ref{example4}).  Then $F_{SO}(0)=l(k+l)>0$, $F_{SO}(z)\to\infty$ as $z\to\infty$, and
$F_{SO}(1)=2k^2-2kl+k-2-l^2+2l<0$ for $l>k$, so $F_{SO}(z)=0$ has two positive solutions.
From the above discussion these solutions are Einstein metrics, which are different from Jensen's Einstein metrics.
Thus the result follows.
\end{proof}

Next, we apply Proposition \ref{main} for $s=2$ and $t=1$ for the
simplectic case $Sp(k_1+k_2+k_3)/Sp(k_3)$. By the
Lagrange method we obtain the following system of equations:

\begin{center}
\begin{equation}
\label{example6}
\begin{array}{lcr}
x_2x_{23}^2\Bigl((k_1+1)x_{12}^2x_{13}^2+k_2x_1^2x_{13}^2+k_3x_1^2x_{12}^2\Bigr)=
x_1x_{13}^2\Bigl((k_2+1)x_{12}^2x_{23}^2+k_1x_2^2x_{23}^2+k_3x_2^2x_{12}^2\Bigr),\\
\\
2x_{13}\Bigl((k_2+1)x_{12}^2x_{23}^2+k_1x_2^2x_{23}^2+k_3x_2^2x_{12}^2\Bigr)=
x_2x_{23}\Bigl(4(k_1+k_2+k_3+1)x_{12}x_{13}x_{23}-\\
(2k_1+1)x_1x_{13}x_{23}-(2k_2+1)x_2x_{13}x_{23}+2k_3x_{12}^3-
2k_3x_{12}x_{13}^2-2k_3x_{12}x_{23}^2\Bigr),\\
\\
x_{13}\Bigl(4(k_1+k_2+k_3+1)x_{12}x_{13}x_{23}-(2k_1+1)x_1x_{13}x_{23}-\\
(2k_2+1)x_2x_{13}x_{23}+2k_3x_{12}^3-2k_3x_{12}x_{13}^2-2k_3x_{12}x_{23}^2\Bigr)=\\
x_{12}\Bigl(4(k_1+k_2+k_3+1)x_{12}x_{13}x_{23}-(2k_1+1)x_1x_{12}x_{23}+
2k_2x_{13}^3-2k_2x_{13}x_{12}^2-2k_2x_{13}x_{23}^2\Bigr),\\
\\
x_{23}\Bigl(4(k_1+k_2+k_3+1)x_{12}x_{13}x_{23}-(2k_1+1)x_1x_{12}x_{23}+
2k_2x_{13}^3-2k_2x_{13}x_{12}^2-2k_2x_{13}x_{23}^2\Bigr)=\\
x_{13}\Bigl(4(k_1+k_2+k_3+1)x_{12}x_{13}x_{23}-(2k_2+1)x_2x_{12}x_{13}+
2k_1x_{23}^3-2k_1x_{23}x_{12}^2-2k_1x_{23}x_{13}^2\Bigr).
\end{array}
\end{equation}
\end{center}

Assume that $x_{13}=x_{23}=z$. Analogously with
$SO(k_1+k_2+k_3)/SO(k_3)$, system 
(\ref{example6})
reduces to the following form:
\begin{center}
\begin{equation}
\label{example7}
\begin{array}{lcr}
x_1=\frac{2(k_1-k_2)x_{12}+(2k_2+1)x_2}{2k_1+1},\\
x_2=\frac{(k_2+1)z^2x_{12}}{(2k_2+k_1+1)z^2+k_3x_{12}^2},\\
z^4x_{12}^3(k_2+1)(2k_2+1)(k_2-k_1)\Bigl((2k_2+2k_1+1)z^2+2k_3x_{12}^2\Bigr)
\Bigl((k_2+k_1)z^2+k_3x_{12}^2\Bigr)^2=0,\\
\\
2(k_3^2+k_1k_3)x_{12}^4-4k_3(k_3+k_2+k_1+1)zx_{12}^3+(2k_1^2+2k_2^2+4k_1k_2+\\
4k_1k_3+10k_2k_3+2k_1+3k_2+6k_3+1)z^2x_{12}^2-(8k_2k_3+4k_1k_3+12k_2k_1+4+4k_3+\\
8k_1+12k_2+8k_2^2+4k_1^2)z^3x_{12}+(2k_1^2+10k_2k_1+6k_1+8k_2+8k_2^2+2)z^4=0.
\end{array}
\end{equation}
\end{center}

From the third equation of (\ref{example7}) we get that a solution
exists only when $k_1=k_2$. Let $k_1=k_2=k$, ($k\geq 3$), $k_3=l$,
$x_1=x_2=x$, $x_{12}=1$, $x_{13}=x_{23}=z$. Then the original system
(\ref{example6}) reduces to:

\begin{center}
\begin{equation}
\label{example8}
\begin{array}{lcr}
2z(x-1)\Bigl(((3k+1)x-k-1)z^2+lx \Bigr)=0,\\
z\Bigl((4kx+2x-8k-4)z^2+(8k+4l+4)z-(2k+1)x-2k-2l \Bigr)=0.
\end{array}
\end{equation}
\end{center}

If $x=1$ we obtain Jensen's solutions, so assume that $x\neq
1$. Then from the first equation of the system (\ref{example8}) we
get $x=\frac{(k+1)z^2}{(3k+1)z^2+l}$. In particular, $0<x<1$ for any
positive $z$. Substituting this expression to the second equation of
(\ref{example8}) we obtain the equation $F_{Sp}(z)=0$, where

\begin{center}
\begin{equation}
\label{example9}
\begin{array}{lcr}
F_{Sp}(z)=2(10k^2+7k+1)z^4-4(6k^2+3kl+5k+l+1)z^3+\\
(8k^2+14kl+5k+6l+1)z^2-4l(2k+l+1)z+2l(k+l).
\end{array}
\end{equation}
\end{center}

If equation $F_{Sp}(z)=0$ has a solution $z>0$, then we get a
new Einstein $Sp(2k+l)\times Sp(k)\times Sp(k)$-invariant
metric. Some results of calculation are presented in Table
2.

\begin{prop}\label{simplSp}
If $l\ge k\ge 1$ then the space $Sp(2k+l)/Sp(l)$ admits at least four
$Sp(2k+l)\times Sp(k)\times Sp(k)$-invariant Einstein metrics.
\end{prop}
\begin{proof}
We consider the polynomial (\ref{example9}).  Then
$F_{Sp}(0)=2l(k+l)>0$, $F_{Sp}(z)\to\infty$ as $z\to\infty$, and
$F_{Sp}(1)=4k^2-4kl-k-2l-1-2l^2<0$ for $l\ge k$, so $F_{Sp}(z)=0$ has
two positive solutions. From the above discussion these solutions are
Einstein metrics, which are different from Jensen's Einstein
metrics. Thus the result follows.
\end{proof}

Note that Propositions \ref{simplSO} and \ref{simplSp} are special cases of Theorems \ref{T1} and \ref{T2} stated
in the Introduction.

\begin{table}[t]
\begin{center}
\begin{tabular}{|c||c|c|c|c|c|c|c|c|c|c|c|c|c|c|c|c|c|c|c|}
\hline
$k$&3&4&5&6&7&8&9&10&11&12&13&14&15&16&17&18&19&20\\
$l$& & & & & & & & & & & & & & & & & & \\
\hline \hline
1 &2&2&2&0&0&0&0&0&0&0&0&0&0&0&0&0&0&0\\
\hline
2 &2&2&2&2&2&0&0&0&0&0&0&0&0&0&0&0&0&0\\
\hline
3 &2&2&2&2&2&2&2&0&0&0&0&0&0&0&0&0&0&0\\
\hline
4 &2&2&2&2&2&2&2&2&2&2&0&0&0&0&0&0&0&0\\
\hline
5 &2&2&2&2&2&2&2&2&2&2&2&2&0&0&0&0&0&0\\
\hline
6 &2&2&2&2&2&2&2&2&2&2&2&2&2&2&2&0&0&0\\
\hline
7 &2&2&2&2&2&2&2&2&2&2&2&2&2&2&2&2&2&0\\
\hline
8 &2&2&2&2&2&2&2&2&2&2&2&2&2&2&2&2&2&2\\
\hline
9 &2&2&2&2&2&2&2&2&2&2&2&2&2&2&2&2&2&2\\
\hline
10&2&2&2&2&2&2&2&2&2&2&2&2&2&2&2&2&2&2\\
\hline
11&2&2&2&2&2&2&2&2&2&2&2&2&2&2&2&2&2&2\\
\hline
12&2&2&2&2&2&2&2&2&2&2&2&2&2&2&2&2&2&2\\
\hline
13&2&2&2&2&2&2&2&2&2&2&2&2&2&2&2&2&2&2\\
\hline
14&2&2&2&2&2&2&2&2&2&2&2&2&2&2&2&2&2&2\\
\hline
15&2&2&2&2&2&2&2&2&2&2&2&2&2&2&2&2&2&2\\
\hline
16&2&2&2&2&2&2&2&2&2&2&2&2&2&2&2&2&2&2\\
\hline
17&2&2&2&2&2&2&2&2&2&2&2&2&2&2&2&2&2&2\\
\hline
18&2&2&2&2&2&2&2&2&2&2&2&2&2&2&2&2&2&2\\
\hline
19&2&2&2&2&2&2&2&2&2&2&2&2&2&2&2&2&2&2\\
\hline
20&2&2&2&2&2&2&2&2&2&2&2&2&2&2&2&2&2&2\\
\hline
\end{tabular}

\vspace {5mm}

{\bf Table 2.} The number of positive solutions of the equation $F_{Sp}(z)=0$
for various $(k,l)$.  (New $Sp(2k+l)\times Sp(k)\times
Sp(k)$-invariant Einstein metrics on $Sp(2k+l)/Sp(l)$).
\end{center}
\end{table}

\section{Another construction for searching for Einstein metrics}\label{Section5}
As shown on Section \ref{Section4} we can find by our method new invariant Einstein metrics on
$SO(k_1\times k_2\times k_3)/SO(k_3)$ (resp. $Sp(k_1\times k_2\times k_3)/Sp(k_3)$), only when $k_1=k_2$.
It would be reasonable to extend this idea to the case $t=1$ and $k_1=k_2=\cdots =k_s=k$ ($s\ge 2$), $k_{s+1}=l$.
Then $n=sk+l$.  If we choose $L'=(SO(k))^s$ (resp. $L'=(Sp(k))^s$), then by Lemma \ref{L2} the
set of $SO(n)\times (SO(k))^s$ (resp. $Sp(n)\times (Sp(k))^s$)-invariant metric
depends on $(s^2+3s)/2$ parameters, which makes the problem difficult for big values $s$.

However, if we choose $L'=N_{SO(sk)}(SO(k))^s$, the normalizer of $(SO(k))^s$ in $SO(sk)$ (resp. $L'=N_{Sp(sk)}(Sp(k))^s$),
(these are extensions of
$(SO(k))^s$ and $(Sp(k))^s$ respectively by a discrete subgroup), then the number of parameters of corresponding
$SO(n)\times L$-invariant metrics reduces to three.  More precisely, the following lemma holds:

\begin{lemma}  If $L'$ is chosen as above, and $K=L'\times H'$, where $H'=SO(l)$ or $Sp(l)$, then
we have a decomposition of $\pp$ into a sum of $\Ad (K)$-invariant and $\Ad (K)$-irreducible submodules
\begin{equation}\label{decomp2}
\pp =\widetilde\pp _1\oplus \widetilde\pp _2\oplus\widetilde\pp _3,
\end{equation}
where $\widetilde\pp _1=\oplus _{i=1}^s\pp _i$,
$\widetilde\pp _2=\oplus _{1\le i<j\le s}\pp _{(i,j)}$, and $\widetilde\pp _3=\oplus _{i=1}^s\pp _{(i,s+1)}$
(cf. (\ref{decomp1})).
The submodules $\widetilde\pp _1, \widetilde\pp _2$ and $\widetilde\pp _3$ are pairwise inequivalent, therefore
any $\Ad (K)$-invariant  inner product of $\pp$ is given by
\begin{equation}\label{last}
(\cdot,\cdot)=
x\cdot \langle \cdot,\cdot \rangle |_{\widetilde\pp _1}+
y\cdot \langle \cdot,\cdot \rangle |_{\widetilde\pp _2}+
z\cdot \langle \cdot,\cdot \rangle |_{\widetilde\pp _3}.
\end{equation}
\end{lemma}
\begin{proof}  For any $1\le i<j\le s$ any two of the submodules $\pp _i$ and $\pp _j$ are interchanged by
$\Ad (a)$, for some $a\in L$.  Similarly, any two of $\pp _{(i,s+1)}$ and $\pp _{(j,s+1)}$ ($1\le i,j\le s$)
are interchanged, and any
two of $\pp _{(i,j)}$ and $\pp _{(i',j')}$ ($1\le i<j\le s$, $1\le i'<j'\le s$).
Therefore decomposition (\ref{decomp2}) follows.
The other statements are obvious.
\end{proof}
Next, we compute the scalar curvature for metric (\ref{last}).
\begin{prop}\label{scalvariant}
The scalar curvature $S$ of an $\Ad (K)$-invariant metric (\ref{last}) has the form

$$
\frac{8(n-2)}{sk}\cdot S=
(k-1)(k-2)\cdot \frac{1}{x}+
(s-1)k((s+2)k-4)\cdot \frac{1}{y}+
4(ks+l-2)l \cdot \frac{1}{z}
$$
\begin{equation}\label{scalvariant1}
-\left(
(s-1)k(k-1) \cdot \frac{x}{y^2}+
(k-1)l \cdot \frac{x}{z^2}+
(s-1)kl \cdot \frac{y}{z^2}
\right)
\end{equation}
with volume condition
$x^\frac{sk(k-1)}{2}y^\frac{s(s-1)k^2}{2}z^{skl}=\ \C$ for the orthogonal case,
and
$$
\frac{4(n+1)}{sk}\cdot S=
(k+1)(2k+1)\cdot \frac{1}{x}+
2(s-1)k((s+2)k+2)\cdot \frac{1}{y}+
8(ks+l+1)l \cdot \frac{1}{z}
$$
\begin{equation}\label{scalvariant2}
-\left(
(s-1)k(2k+1) \cdot \frac{x}{y^2}+
(2k+1)l \cdot \frac{x}{z^2}+
2(s-1)kl \cdot \frac{y}{z^2}
\right)
\end{equation}
with volume condition $x^{s(2k+1)k}y^{2s(s-1)k^2}z^{4skl}=\ \C$
for the symplectic case.
\end{prop}
\begin{proof}
Metric (\ref{last}) is a special case of metric (\ref{1}) for which
the scalar curvature was obtained in Proposition \ref{main}. We
apply these expressions for $t=1$, $k_1=\cdots =k_s=k$, $k_{s+1}=l$,
and $x_a=x$ ($1\le a\le s$), $x_{a,b}=y$ ($1\le a<b\le s$),
$x_{a,s+1}=z$ ($1\le a\le s$) to obtain
$$
\sum\limits_{a=1}^s\frac{k_a(k_a-1)(k_a-2)}{8(n-2)}\cdot \frac{1}{x_a}=
\frac{sk(k-1)(k-2)}{8(n-2)}\cdot \frac{1}{x};
$$
$$
\sum\limits_{1\leq a<b\leq s+t} \frac{k_ak_b}{x_{(a,b)}}=
\sum\limits_{1\leq a<b\leq s} \frac{k_ak_b}{x_{(a,b)}}+
\sum\limits_{1\leq a\leq s, b=s+1} \frac{k_ak_b}{x_{(a,b)}}=
\frac{s(s-1)k^2}{2}\cdot \frac{1}{y}+
skl \cdot \frac{1}{z};
$$
$$
\sum\limits_{1\leq a\leq s, \,a+1\leq b\leq s+t} k_ak_b(k_a-1)
\frac{x_a}{x_{(a,b)}^2}=
\sum\limits_{1\leq a<b\leq s} k_ak_b(k_a-1)
\frac{x_a}{x_{(a,b)}^2}+
\sum\limits_{1\leq a\leq s, b=s+1} k_ak_b(k_a-1)
\frac{x_a}{x_{(a,b)}^2}=
$$
$$
\frac{s(s-1)k^2(k-1)}{2} \cdot \frac{x}{y^2}+
sk(k-1)l \cdot \frac{x}{z^2};
$$
$$
\sum\limits_{1\leq a<b\leq s} k_ak_b(k_b-1)
\frac{x_b}{x_{(a,b)}^2}=
\frac{s(s-1)k^2(k-1)}{2}\cdot \frac{x}{y^2};
$$
$$
\sum\limits_{1\leq a <b<c \leq s+t} k_ak_bk_c
\left(
\frac{x_{(a,b)}}{x_{(a,c)}x_{(b,c)}}+
\frac{x_{(a,c)}}{x_{(a,b)}x_{(b,c)}}+
\frac{x_{(b,c)}}{x_{(a,b)}x_{(a,c)}}
\right)    =
$$
$$
\sum\limits_{1\leq a <b<c \leq s} k_ak_bk_c
\left(
\frac{x_{(a,b)}}{x_{(a,c)}x_{(b,c)}}+
\frac{x_{(a,c)}}{x_{(a,b)}x_{(b,c)}}+
\frac{x_{(b,c)}}{x_{(a,b)}x_{(a,c)}}
\right)
$$
$$
+\ \sum\limits_{1\leq a <b\leq s, c=s+1} k_ak_bk_c
\left(
\frac{x_{(a,b)}}{x_{(a,c)}x_{(b,c)}}+
\frac{x_{(a,c)}}{x_{(a,b)}x_{(b,c)}}+
\frac{x_{(b,c)}}{x_{(a,b)}x_{(a,c)}}
\right)    =
$$
$$
\frac{s(s-1)(s-2)k^3}{2}\cdot \frac{1}{y}+
\frac{s(s-1)k^2l}{2}\cdot \left( \frac{y}{z^2}+\frac{2}{y} \right),
$$
for the orthogonal case, so equation (\ref{scalvariant1}) is obtained.
For the symplectic case we have that
$$
S=\sum\limits_{a=1}^s\frac{k_a(k_a+1)(2k_a+1)}{4(n+1)}\cdot \frac{1}{x_a}+
2
\sum\limits_{1\leq a<b\leq s+t} \frac{k_ak_b}{x_{(a,b)}}
$$
$$
-\frac{1}{4(n+1)}
\sum\limits_{1\leq a\leq s,\, a+1\leq b\leq s+t} k_ak_b(2k_a+1)
\frac{x_a}{x_{(a,b)}^2}-
\frac{1}{4(n+1)}
\sum\limits_{1\leq a<b \leq s} k_ak_b(2k_b+1)
\frac{x_b}{x_{(a,b)}^2}
$$
$$
- \frac{1}{n+1}
\sum\limits_{1\leq a <b<c \leq s+t} k_ak_bk_c
\left(
\frac{x_{(a,b)}}{x_{(a,c)}x_{(b,c)}}+
\frac{x_{(a,c)}}{x_{(a,b)}x_{(b,c)}}+
\frac{x_{(b,c)}}{x_{(a,b)}x_{(a,c)}}
\right),
$$
which implies equation (\ref{scalvariant2}), analogously to the orthogonal case.
The dimensions of $\widetilde\pp _1, \widetilde\pp _2$ and $\widetilde\pp _3$ are
for the orthogonal case $\frac{sk(k-1)}2, \frac{s(s-1)k^2}2$ and $skl$ respectively,
and for the symplectic case $sk(2k+1), 2s(s-1)k^2$ and $4skl$ respectively.
Thus the volume conditions are obtained, and the proof is completed.
\end{proof}

\medskip
In order to find the critical points of the scalar curvature $S$ for the above two cases,
note that
$\frac{8(n-2)}{sk}\cdot S$ (respectively,$\frac{4(n+1)}{sk}\cdot S$)
and the volume are  of the form
\begin{equation}\label{vs1}
F(x,y,z)=\frac{a}{x}+\frac{b}{y}+\frac{c}{z}-d\frac{x}{y^2}-e\frac{x}{z^2}-f\frac{y}{z^2},
\quad G(x,y,z)=x^py^qz^r,
\end{equation}
where the constant $a$, $b$, $c$, $d$, $e$, $f$, $p$, $q$ and $r$
are positive, and
\begin{equation}\label{vs2}
d=\frac{pb-qa}{q+2p}, \quad f=\frac{qe}{p}.
\end{equation}
We need to consider the following problem: Find all the critical
points (with positive coordinates) of $F(x,y,z)$ under the
constraint $G(x,y,z)=\C$. This is a Lagrange-type
problem.

\begin{lemma}\label{vs1.1}
The critical points of the function $F(x,y,z)$ with positive
$x,y,z$ under the restriction $G(x,y,z)=\C$ satisfy the following
equations:

1) If $x=y$, then
$$
r(a+d)z^2-pcxz+e(2p+2q+r)x^2=0;
$$

2) If $x \neq y$, then
$$
x=\frac{aqyz^2}{pfy^2+d(q+2p)z^2},
$$
$$
(2d(q+2p)+bq)drz^4-(q+2p)cdqyz^3+ (2d(r+q)(q+2p)+(r+2p)aq)fy^2z^2-
$$
$$
cfpqy^3z+(r+2q)f^2py^4=0.
$$
If in addition $d(q+2p)>aq$, then $y>x$.
\end{lemma}
\begin{proof}
It is easy to see that the problem reduces to the following
system:
\begin{equation}\label{ur1} \left\{
\begin{array}{lcr}
q\big(-\frac{a}{x}-d\frac{x}{y^2}-e\frac{x}{z^2}\big)=p\big(-\frac{b}{y}+2d\frac{x}{y^2}-f\frac{y}{z^2}\big),\\
\\
r\big(-\frac{a}{x}-d\frac{x}{y^2}-e\frac{x}{z^2}\big)=p\big(-\frac{c}{z}+2e\frac{x}{z^2}+2f\frac{y}{z^2}\big).
\end{array}
\right.
\end{equation}
From the first equation of  (\ref{ur1}) we get:
$$
pf\frac{y-x}{z^2}=\frac{aq(y-x)(y-\frac{d(q+2p)}{aq}x)}{xy^2}.
$$

If $x=y$ we easily obtain that
$$
r(a+d)z^2-pcxz+e(2p+2q+r)x^2=0.
$$
If $x\neq y$ we obtain
$$
x=\frac{aqyz^2}{pfy^2+d(q+2p)z^2},
$$
which implies that $x>0$ for any $z>0, y>0$. If $x\neq y$ then
$$
\frac{y}{x}>\frac{d(q+2p)}{aq},
$$
therefore, if  $d(q+2p)>aq$, then $y>x$.

Substituting the above expression for $x$ in the second equation of
(\ref{ur1}), we obtain the Einstein equation
$$
(2d(q+2p)+bq)drz^4-(q+2p)cdqyz^3+ (2d(r+q)(q+2p)+(r+2p)aq)fy^2z^2-
$$
$$
cfpqy^3z+(r+2q)f^2py^4=0.
$$
\end{proof}

Let
$$
P(u)=(2d(q+2p)+bq)dru^4-(q+2p)cdqu^3+ (2d(r+q)(q+2p)+(r+2p)aq)fu^2-
$$
\begin{equation}\label{P}
cfpqu+(r+2q)f^2p.
\end{equation}
It is clear that the equation $P(z/y)=0$ is equivalent to the last one of
Lemma \ref{vs1.1}. The following simple lemma will be
used for the proof of the main theorems

\begin{lemma}\label{vs2.1} If $P(1)<0$ then  equation $P(u)=0$ has at least two
positive solutions.
\end{lemma}
\begin{proof}
It is evident from the facts that
$P(0)=(r+2q)f^2p>0$ and  $P(u)\to +\infty \mbox{ when }u\to +\infty.$
\end{proof}

\section{Proof of the main results}\label{Section6}

\begin{proof}[Proof of Theorem \ref{T1}]
We apply Lemma \ref{vs1.1}  for values of $a, b,
c, d, e, f, p, q$ and $r$ taken from the orthogonal case  of  Proposition
\ref{scalvariant}. If $x=y$ then one gets obviously Jensen's
discussed in Section \ref{SectionJen}. If $x\ne y$ then
for the polynomial (\ref{P}) it is that
$$
P(1)=\frac{1}{2}s^2k^4l(k-1)(s-1)^2(sk^2-skl+k-2-l^2+2l).
$$
It it easy to check that $sk^2-skl+k-2-l^2+2l<0$ for $l>k$, $k,s\ge
2$, thus $P(1)<0$. By Lemma \ref{vs2.1} the equation $P(u)=0$ has
at least two positive solutions so we obtain at least two new
invariant Einstein metrics.  Since in this case $d(q+2p)>aq$, then $y>x$ for these
new metrics.
\end{proof}

\medskip
\begin{proof}[Proof of Theorem \ref{T2}]
We use Lemma \ref{vs1.1} for values of $a, b,
c, d, e, f, p, q$ and $r$ taken from the symplectic case  of  Proposition
\ref{scalvariant}. If $x=y$ then we obtain the two Jensen's
solutions as derived in Section \ref{SectionJen}. If $x\ne y$ then
for the polynomial (\ref{P}) it is that
$$
P(1)=8s^2k^4l(2k+1)(s-1)^2(2sk^2-2skl-k-2l-1-2l^2).
$$
It it easy to check that $2sk^2-2skl-k-2l-1-2l^2<0$ for $l\ge k$,
$s\ge 2$, thus $P(1)<0$. By Lemma \ref{vs2.1} the equation
$P(u)=0$ has at least two positive solutions, so obtain at
least two new invariant Einstein metrics.  Since in this case it is
also $d(q+2p)>aq$, then $y>x$ for these
new metrics.
\end{proof}

\medskip
 \begin{proof}[Proof of Theorem \ref{T3}]
 Fix a positive integer $p$ and choose positive integers $n, l$ such that $n-l$ has at least $p$ different prime factors
 $a_1, a_2, \dots, a_p$ with $a_i<l$ ($i=1,\dots ,p$).  Take $k$ any of the $a_i$'s, and positive integer
 $s$ so that $n-l=sk$.
 For this choice of $k, l, s$ we use Theorems  \ref{T1} and \ref{T2}, and obtain that the homogenous spaces
 $SO(n)/SO(l)$ and $Sp(n)/Sp(l)$ admit at least two $SO(n)\times (SO(k))^s$ (resp. $Sp(n)\times (Sp(k))^s$)-invariant
 Einstein metrics which are not invariant under the group $SO(n)\times SO(n-l)$ (resp. $Sp(n)\times Sp(n-l)$), that is they are not
 Jensen's metrics.
 It is easy to see that for different choices of $k=a_i$ we obtain pairwise different metrics (because they have different
 full motion groups).
 Therefore, we obtain at least $2p$ pairwise different $SO(n)$ (resp. $Sp(n)$)-invariant Einstein metrics on the
 Stiefel manifold $SO(n)/SO(l)$ (resp. the space $Sp(n)/Sp(l)$).
\end{proof}

\medskip
It would be an interesting problem to investigate the nature of the invariant Einstein metrics of given volume
on the spaces $SO(n)/SO(l)$ and $Sp(n)/Sp(l)$ (cf. Theorems \ref{T1} and \ref{T2}),
as critical points of the scalar curvature functional curvature
$S$. For instance, by analysing the Hessian of $S$ at the critical points.
Of course, this would require having explicit solutions of the algebraic systems
of equations obtained from the Einstein equation.
Another interesting problem is to find metrics with maximal and minimal values
of the scalar curvature $S$ among all invariant Einstein metrics of fixed volume
on the spaces $SO(n)/SO(l)$ and $Sp(n)/Sp(l)$.

\end{document}